\documentclass{article}
\usepackage{amsmath,amssymb}
\usepackage{graphicx}
\usepackage{xypic}

\newtheorem{prop}{Proposition} 
\newtheorem{theo}[prop]{Theorem}  
\newtheorem{lemm}[prop]{Lemma}    
\newtheorem{coro}[prop]{Corollary}    

\def\proof{\noindent\textit{Proof. }}
\newcounter{example}
\def\example{\addtocounter{example}{1}\noindent\textit{Example \arabic{example}. }}

\newcommand{\pic}[4]{\vspace{1ex}\setlength{\unitlength}{1cm}
\begin{picture}(0,#3)(5,.5)
\put(#2){\includegraphics[#4]{#1.eps}}
\end{picture}\vspace{1ex}}

\def\C{{\mathbb C}}

\def\R{{\mathbb R}}

\def\sgn{\,\!\mbox{sgn}\,}

\def\DS{\displaystyle}
\def\qed{\ \hfill\raise3.5pt\hbox{\framebox[1ex]{\ }}}
\def\vd{\vec\delta}
\def\vj{\vec\jmath\hspace{.3ex}}
\def\vn{\vec n}

\newcommand{\X}[1]{X^{(#1)}}
\newcommand{\Xt}[1]{\widetilde{X}^{(#1)}}
 
\newcommand{\stX}[1]{{\overset{*}{X}}^{\raisebox{-1.6ex}{\scriptsize$(#1)$}}}
\newcommand{\stXt}[1]{\overset{\ *}{\widetilde{X}}^{\raisebox{-2.3ex}{\scriptsize$(#1)$}}}
\newcommand{\ststX}[1]{{\overset{**}{X}}^{\raisebox{-1.6ex}{\scriptsize$(#1)$}}}

\begin{document}

\begin{center}\Large
\textbf{On Sturm-Liouville Equations with\\ Several Spectral Parameters}\\
\large R. Michael Porter \\   \today
 \end{center} 

\noindent\textbf{Abstract.}  We give explicit formulas for a pair of
linearly independent solutions of
$(py')'(x)+q(x)=(\lambda_1r_1(x)+\cdots+\lambda_dr_d(x))y(x)$, thus
generalizing to arbitrary $d$ previously known formulas for $d=1$.
These are power series in the spectral parameters
$\lambda_1,\dots,\lambda_d$ (real or complex), with coefficients which
are functions on the interval of definition of the differential
equation.  The coefficients are obtained recursively using indefinite
integrals involving the coefficients of lower degree.  Examples are
provided in which these formulas are used to solve numerically some
boundary value problems for $d=2$, as well as an application to
transmission and reflectance in optics.



\setcounter{section}{-1} 
\section{Introduction}
We will consider the second-order linear differential equation
\begin{equation}  \label{eq:severaleigen}
  (py')' + qy = (\lambda_1 r_1 + \lambda_2 r_2 +\cdots+ \lambda_d r_d )y
\end{equation}
on a real interval $x_1\le x\le x_2$, where $p,q,r_1,\dots,r_d$ are given
functions and $\lambda_1,\dots,\lambda_d$ are unknown parameters.  Let
some appropriate boundary conditions be imposed at $x_1$ and $x_2$.  Then a
\emph{spectral problem} consists of determining the subset
$(\lambda_1,\dots,\lambda_d)\subset\R^d$ (or $\C^d$) for which there
exists a solution $y$ of  (\ref{eq:severaleigen}) which satisfies those
boundary conditions. While there is a vast literature on spectral
theory for general differential equations and on numerical methods
specifically developed for (\ref{eq:severaleigen}) for $d=1$---indeed, the
expression ``spectral problem'' commonly implies a single
$\lambda_1$---   there is considerably less available concerning
several parameters.  One may find qualitative results on this subject
in \cite{At1,AtMi,BrSl,Sl,Tu,Vol}.  For $d=2$ some properties of the
eigencurves are set forth in \cite[chapter 6]{AtMi}.

An approach for solving spectral problems for $d=1$ was presented in
\cite{KrCV08,KrP2010} which produces two explicit power series in the
variable $\lambda_1$ with coefficients which are functions on
$[x_1,\,x_2]$.  These series represent two functions $y=u_1(x)$,
$y=u_2(x)$ parametrized by $\lambda_1$ which are linearly independent
solutions of (\ref{eq:severaleigen}). There are similar power series
for the derivatives $u'_1(x)$, $u'_2(x)$.  By evaluating these series
with $x$ at the endpoints $x_1$, $x_2$ we obtain power series in
$\lambda_1$ for the boundary values, which upon substitution
in the boundary conditions produce a ``characteristic function'' whose
zeroes are the eigenvalues of the spectral problem. (It is not
necessary for the boundary conditions to be linear for this procedure
to apply.) These series representations have applications beyond
spectral problems; for example they provide an effective method for
solving initial value problems.

Since its appearance in 2008, consequences of this SPPS (spectral
parameter power series) representation have been investigated in many
directions. These include completeness properties of the ``formal
powers'' used to define the coefficients of the power series
\cite{KrCMA2011,KMoT}; relationship to transmutation operators,
Darboux and other transformations, and Goursat problems
\cite{KhKrTT,KrT2012,KrT2015hyp,KrT2015trans}; extension to other number
systems (quaternions, etc) \cite{CKr1,CKr2,KrT2015hyp} and equations of
higher order \cite{KKB}; relaxation of regularity conditions on the
coefficients of the differential equation \cite{BCKh,CKrT}.  Further, there
have appeared numerous applications to problems in physics and
engineering \cite{CKKO,KhKrR,KhS,KhT,KrV2011} as well as in complex analysis
\cite{BrP2011,KRP2011}.

In dealing with physics or engineering models which involve a
Sturm-Liouville problem containing several eigenvalues $\lambda_i$, it
is common practice to fix all but one of them, and then solve the
spectral problem for the remaining one. This appears to be due to the
difficulties of existing methods of handling more than one spectral
parameter. In this paper we work out the SPPS coefficients
corresponding to (\ref{eq:severaleigen}) for arbitrary $d\ge1$. 
This permits treating the spectral parameters in unified way. We
give some numerical examples with $d=2$, and then an application to a
problem of transmittance of an electromagnetic wave through an
inhomogeneous layer, in which the two parameters correspond to
physical characteristics of the phenomenon.

\section{Formal powers \label{sec:formal}}

 The Sturm-Liouville linear
differential expression on the left-hand side of
(\ref{eq:severaleigen}) will be denoted by
\begin{equation}  \label{eq:defL}
  Ly = (py')' + qy.
\end{equation}
Throughout this paper $p,q,r_1,\dots,r_d$ will denote real or complex
valued functions on the closed interval $[x_1,\,x_2]$.  In this section
we are interested in describing the procedure for constructing the
SPPS representation of solutions, while questions of convergence and
regularity will be deferred to the next section. A basepoint $x_0$ is
fixed in  $[x_1,\,x_2]$. For convenience we will use the notation
\[  g=\int\!f \]
to mean
\[    g(x)=\int_{x_0}^xf(s)\,ds
\]
for any function $f$ under consideration, inasmuch as we will have no
use for other limits of integration.  In the following we will set up
a notation for describing sums of finitely nested integrals of the
form
\begin{eqnarray}
  \cdots   \int r_{i_n}u_0^2  \int \frac{1}{pu_0^2}\cdots 
           \int r_{i_2}u_0^2  \int \frac{1}{pu_0^2}  \int r_{i_1}u_0^2  
   \label{eq:Xtform},  && \\
   \cdots  \int r_{i_n}u_0^2  \int \frac{1}{pu_0^2}\cdots  
           \int r_{i_2}u_0^2  \int \frac{1}{pu_0^2}  \int r_{i_1}u_0^2
           \int \frac{1}{pu_0^2} .  \label{eq:Xform}  
\end{eqnarray} 

\subsection{Construction of $\Xt{\vj}$\label{subsec:u1}}

For simplicity of handling the indices, we will begin with the form
(\ref{eq:Xtform}) which produces the family of functions called
$\Xt{\vj}$ (the notation follows that generally used in the SPPS
literature). Here $\vj=(j_1,\dots,j_d)$ is a multiindex with integral
entries.  It has $d$ predecessors given by
\[  \ \vj - \vd_i = (j_1,\dots, j_{i-1},j_i-1,j_{i+1},\dots,j_d)
\]
where $\vd_i=(0,\dots,0,1,0,\dots,0)$ is the $i$-th standard basis
vector.  We will say that $\vj$ is an \emph{admissible} multiindex
when at most one of its entries $j_i$ is odd.  An admissible $\vj$ is
called \emph{even} or \emph{odd} according to the parity of 
\[  |\vj|=\sum_{i=1}^d j_i; \]
i.e., it is odd when exactly one $j_i$ is odd.
We start from the constant function
\begin{equation}   \label{eq:defXt0}
   \Xt{\vec0}(x)=1
\end{equation}
for all $x$, where $\vec0=(0,\dots,0)$. For definiteness we set
$\Xt{\vj}(x)=0$ whenever $j_i<0$ for some $i$.  Then we define the
formal power $\Xt{\vj}$ for admissible $\vj$ with 
nonnegative indices   in the following recursive manner:
\begin{equation}   \label{eq:defXt}
 \Xt{\vj} = \left\{ \begin{array}{ll}
 \DS |\vj| \int r_iu_0^2 \,\Xt{\vj-\vd_i},\quad  &
      \vj\mbox{ odd},  \\[2ex]
 \DS |\vj| \int \frac{1}{pu_0^2}\sum_{i=1}^d\Xt{\vj-\vd_i}, \quad&
      \vj\mbox{ even}.\\
 \end{array} \right.  
\end{equation}

\begin{figure}[!b]  \centering  
\newcommand{\rowee}[1]{ \Xt{#1} \ar[r]^{r_2}\ar[d]_{r_1} }
\newcommand{\roweo}[1]{ \Xt{#1} \ar[r]^{p} }
\newcommand{\rowoe}[1]{ \Xt{#1} \ar[d]_{p} }
\[ \xymatrix{
  \rowee{0,0} & \roweo{0,1}  & \rowee{0,2}    & \roweo{0,3} & \rowee{0,4} \ar[d]_{r_1}  & \Xt{0,5}\\ 
   \rowoe{1,0} & & \rowoe{1,2} & & \rowoe{1,4} \\
  \rowee{2,0}\ar[d]_{r_1}  & \roweo{2,1}  & \rowee{2,2}  & \roweo{2,3} & \rowee{2,4}   & \Xt{2,5}\\ 
  \rowoe{3,0} & & \rowoe{3,2} & & \rowoe{3,4} \\
  \rowee{4,0} \ar[d]_{r_1} & \roweo{4,1}  & \rowee{4,2} \ar[d]_{r_1} & \roweo{4,3} & \rowee{4,4} & \Xt{4,5}\\ 
  \Xt{5,0} & & \Xt{5,2} & & \Xt{5,4} \\
 }
\]
  \caption{Construction of $\Xt{\vj}$.}
  \label{fig:Xt}
\end{figure}

  These are all functions on $[x_1,x_2]$. Note
that when $\vj$ is odd, the index $i$ referred to in the first clause
of (\ref{eq:defXt})  is unambiguously defined.  The interdependencies of the
$\Xt{\vj}$ are illustrated for $d=2$ in Figure \ref{fig:Xt}. We will
say that the \emph{degree} of $\Xt{\vj}$ is $|j|$.

\begin{samepage}
To motivate to some extent what we have done we give the following relationship.
\begin{lemm} \label{lemm:Lu0Xt}
  Let $u_0$ be a nonvanishing function on $[x_1,x_2]$ and suppose that
  $Lu_0=0$. Then for any nonnegative even multiindex $2\vn$,
 \[  L(u_0\Xt{2\vn}) = 
   2|\vn|(2|\vn|-1)\, u_0 \sum_{i=1}^d
               r_i \Xt{2\vn -2\vd_i}.
\]
\end{lemm}
\end{samepage}
\proof As a consequence of $Lu_0=0$, is easily seen that the operator
$L$ admits the Polya factorization \cite{KelPet}
\[  L = \frac{1}{u_0}\, \partial\, pu_0^2\, \partial\, \frac{1}{u_0}
\]
where $\partial=\partial/\partial x$ and the functions in this expression
refer to the corresponding multiplication operators. Thus by (\ref{eq:defXt}),
\begin{eqnarray*}
 L(u_0\Xt{2\vn }) &=& 
  \frac{1}{u_0}\partial\, pu_0^2\,\partial\frac{1}{u_0}\left(
  2|\vn|u_0\int \frac{1}{pu_0^2} \sum_i \Xt{\vec 2n-\vd_i} \right)\\
  &=&  2|\vn| \frac{1}{u_0}\partial \left(
   \sum_i \Xt{2\vn  - \vd_i} )  \right). 
\end{eqnarray*}
A second application of (\ref{eq:defXt}) yields that this is equal to
\[  2|\vn| \frac{1}{u_0} \left(
   (2|\vn|-1)   \sum_i   r_iu_0^2 \Xt{2\vn  - 2\vd_i}
     \right) 
\]
as required. \qed

The number $\tilde c_{\vj}$ of summands of the form (\ref{eq:Xtform})
comprising $\Xt{\vj}$ is the same as the number of paths which advance
(i.e.\ from predecessors to successors) from $\Xt{\vec0}$ to $\Xt{\vj}$
and can be described as follows. We define $\tilde c_{\vj}=0$ when
some $j_i<0$. Clearly $\tilde c_{\vec0}=1$. Then by (\ref{eq:defXt})
we have recursively
\begin{equation}  \label{eq:ctilde}
 \tilde c_{\vj} = \left\{ 
    \begin{array}{ll}
     \tilde  c_{\vj-\vd_i}, \quad& \vj\mbox{ odd},\\[2ex]
    \DS \sum_{i=1}^d  \tilde c_{\vj-\vd_i} , \quad& \vj\mbox{ even}.
    \end{array} \right.
\end{equation} 
By induction via predecessors it is readily seen that
\begin{equation}  \label{eq:ctildebound}
 \tilde c_{\vj} = \frac{ [ \frac{|\vj|}{2} ]! }
     { [\frac{j_1}{2}]!  [\frac{j_2}{2}]! \cdots  [\frac{j_d}{2}]!} .
\end{equation}

Consider a single nested integral appearing as a summand in $\Xt{\vj}$. The number
of integrations following division by $pu_0^2$ is $[|\vj|/2]$, while
for each $i$ ($1\le i\le d$), the number of integrations which follow
a multiplication by $r_iu_0^2$ is easily seen to be $[(j_i+1)/2]$.
Here and always $[a]$ means the least integer no greater than the real
number $a$.  One may verify that these indeed sum to $|\vj|$. 
Define
\begin{equation}  \label{eq:M0Mi}
 M_0 = \sup_{[x_1,\,x_2]} \frac{1}{|pu_0^2|},\quad  M_i=\sup_{[x_1,\,x_2]} |r_iu_0^2| . 
\end{equation}
\begin{lemm}  
The formal powers $\Xt{\vj}$ satisfy the growth condition
\begin{equation}  \label{eq:Xtbound}
  |\Xt{\vj}(x)| \le \tilde c_{\vj} \,
         M_0^{\left[\frac{|\vj|}{2}\right]}
         M_1^{\left[\frac{j_1+1}{2}\right]}
         M_2^{\left[\frac{j_2+1}{2}\right]}\cdots
         M_d^{\left[\frac{j_d+1}{2}\right]}  |x-x_0|^{|\vj|}
\end{equation}
for $x_1\le x\le x_2$.
\end{lemm} 

\proof
Suppose that $\vj$ is even. Then by the inductive hypothesis
\begin{eqnarray*}
   |\Xt{\vj}(x)| &\le&  |\vj| \int_{x_0}^x (\sup \frac{1}{|pu_0^2|}) 
      \sum_{i=1}^d  |\Xt{\vj-\vd_i}(t)| \,dt \\
 &\le& \int_{x_0}^x M_0\left(\sum_{i=1}^d \tilde c_{\vj-\vd_i} 
       M_0^{\left[\frac{|\vj|-1}{2}\right]}
       M_1^{\left[\frac{j_1+1}{2}\right]} \cdots M_i^{\left[\frac{j_i}{2}\right]}   \cdots   M_d^{\left[\frac{j_d+1}{2}\right]}    \right)|\vj| |t-x_0|^{|\vj|-1} \,dt.
\end{eqnarray*}
We then integrate and note that $[j_i/2]=[(j_i+1)/2]$ since all $j_i$ are even,
obtaining the bound
\[  \left( \sum_{i=1}^d \tilde c_{\vj-\vd_i}  
          M_0^{\left[\frac{|\vj|}{2}\right]} M_1^{\left[\frac{j_1+1}{2}\right]}  M_2^{\left[\frac{j_2+1}{2}\right]}
   \cdots  M_d^{\left[\frac{j_d+1}{2}\right]}  
            \right)|x-x_0|^{j+k}
\]
which by (\ref{eq:ctilde}) reduces to (\ref{eq:Xtbound}). The
verification for $|\vj|$ odd is similar and indeed simpler.  \qed

\subsection{Construction of $\X{\vj}$\label{subsec:defX}}

The construction of $\X{\vj}$ in terms of nested integrals of the form
(\ref{eq:Xform}) is analogous to that of $\Xt{\vj}$. However, there
are some notational complications. The indices could be handled in
various ways; our choice, perhaps purist, is as follows. Now the $\vj$
will have entries with common fractional part $j_i-[j_i]=1/d$.  We
will call $\vj$ admissible when at most one of the the integral parts
$[j_i]$ is odd, while the parity of $\vj$ is again that of the integer $|\vj|$.  

To start the recursion we use the predecessors of $(1/d)\vec1$ , i.e.
$\vj=(1/d,1/d,\dots,1/d,-1+1/d,1/d,\dots,1/d)$, defining  the constant functions 
\newcommand{\Xh}[2]{\X{#1\frac{1}{2},#2\frac{1}{2}}}  
 \begin{equation}   \label{eq:defX0}
  \X{ (1/d)\vec{1}-\vd_i}(x) = \frac{1}{d}
\end{equation}
for $i=1,\dots,d$.  These are functions
of degree 0. For definiteness we define $\X{\vj}(x)=0$ whenever some
$j_i<0$ except as specified by (\ref{eq:defX0}).  The formal powers
$\X{\vj}$ for the remaining admissible $\vj\ge0$ are defined by
\begin{equation}\label{eq:defX}
 \X{\vj} = \left\{ \begin{array}{ll}
 \DS  |\vj| \int r_iu_0^2 \, \X{\vj-\delta_i},\quad& \vj\mbox{ even},\\[2ex]
 \DS |\vj|  \int \frac{1}{pu_0^2}\, \sum_{i=1}^d  \X{\vj-\delta_i},\quad& \vj\mbox{ odd},
 \end{array} \right.
\end{equation}
as outlined in Figure \ref{fig:X}.  This formula differs from (\ref{eq:defXt}) not
only in the exchange of even and odd, but also in that the indices
and coefficients have different interpretations. One justification for this
notation is the role of the degree $|\vj|$, cf.\ Lemma
\ref{lemm:Xbound}.

\begin{samepage}
Analogously to Lemma \ref{lemm:Lu0Xt} we find
\begin{lemm} \label{lemm:Lu0X}
\[ L(u_0\X{2\vn +\frac{1}{d}\vec1}) =   
(2(|\vn|+1)(2|\vn|)\, u_0 \sum_{i=1}^d
               r_i \,\Xt{2\vn -2\vd_i+\frac{1}{d}\vec1}. 
\]
\end{lemm} 
\end{samepage}
In verifying this it is useful to note that $|(1/d)\vec1|=1$ and to use the linearity of the degree
operator $|\cdot|$.

\begin{figure}[!t]   \centering
  \newcommand{\rowee}[2]{ \Xh{#1}{#2} \ar[r]^{r_2}\ar[d]_{r_1} }
  \newcommand{\roweo}[2]{ \Xh{#1}{#2} \ar[r]^{p} }
  \newcommand{\rowoe}[2]{ \Xh{#1}{#2} \ar[d]_{p} }
\begin{small}
\[ \xymatrix{
  & \Xh{-}{}\ar[d]_{p}\\ 
  \Xh{}{-}\ar[r]_{p} & 
  \rowee{}{}  & \roweo{}{1}  & \rowee{}{2}  & \roweo{}{3} & \Xh{}{4}\ar[d]_{r_1} \\
 & \rowoe{1}{} & & \rowoe{1}{2} & & \rowoe{1}{4} \\
 & \rowee{2}{} & \roweo{2}{1}  & \rowee{2}{2}  & \roweo{2}{3} &  \Xh{4}{4}\ar[d]_{r_1} \\
 & \rowoe{3}{} & & \rowoe{3}{2} & & \rowoe{3}{4} \\
 & \rowee{4}{} & \roweo{4}{1}  & \rowee{4}{2}  & \roweo{4}{3} &  \Xh{4}{4}\ar[d]_{r_1} \\
 & \Xh{5}{} & & \Xh{5}{2} & & \Xh{5}{4} \\
 }
\]  
\end{small}
\caption{Formal powers $\X{\vj}$ for $d=2$.\label{fig:X}}
\end{figure}

The number $c_{\vj}$ of terms in $\X{\vj}$ is determined recursively
by setting $c_{(1/d)\vec1-\vd}=1$, while otherwise $c_{\vj}=0$ if some
$j_i<0$, and then
\begin{equation}  \label{eq:c}
 \tilde c_{\vj} = \left\{ 
    \begin{array}{ll}
     \tilde  c_{\vj-\vd_i}, \quad& \vj\mbox{ even},\\[2ex]
    \DS \sum_{i=1}^d  \tilde c_{\vj-\vd_i} , \quad& \vj\mbox{ odd},
    \end{array} \right.
\end{equation}
which is analogous (\ref{eq:ctilde}) but again with a different
interpretation of the indices.  The number of integrations following
division by $pu_0^2$ is now $[(|\vj|+1)/2]$, and those following
multiplication by $r_iu_0^2$ number $[(j_i-1/d+1)/2]$.  This gives the
growth estimate:
\begin{lemm} \label{lemm:Xbound}
\[  |\X{\vj}(x)| \le  c_{\vj} M_0^{\left[\frac{\vj+1}{2}\right]}
                         M_1^{\left[\frac{j_1-1/d+1}{2}\right]} \cdots
                         M_d^{\left[\frac{j_d-1/d+1}{2}\right]} |x-x_0|^{|\vj|}.
 \]
\end{lemm}

\section{SPPS series and characteristic function}

\subsection{General solution}
We define the SPPS functions $u_1$, $u_2$ as
\begin{eqnarray}     
  u_1 &=& u_0 \sum_{\vn \ge0}  \frac{1}{(2|\vn|)!}\Xt{2 \vn }
             \lambda_1^{n_1} \cdots\lambda_d^{n_d}, \nonumber  \\
  u_2 &=& u_0\sum_{\vn \ge0}  \frac{1}{(2|\vn|+1)!}\X{2\vn +\frac{1}{d}\vec1}
    \lambda_1^{n_1}\cdots\lambda_d^{n_i}.   \label{eq:defu1u2}
\end{eqnarray}
where the sums are over all nonnegative multiindices $\vn$.  Note that
the degree of $\X{2\vn +\frac{1}{d}\vec1}$ is $2|\vn|+1$. The main
result, which mirrors that of \cite{KrP2010}, is as follows.

\begin{theo} \label{theo:SPPS} Let $p,q,r_1,\dots,r_d$ be continuous
  complex-valued functions of the real variable $x\in[x_0,\,x_1]$,
  with $p$ continuously differentiable and $p(x)\not=0$. Let the
  differential operator $L$ be defined by (\ref{eq:defL}). Then for
  every $\vec\lambda=(\lambda_1,\dots,\lambda_d)\in\C^d$ the two
  series in (\ref{eq:defu1u2}) converge uniformly on
  $x\in[x_0,\,x_1]$, and the functions $u_1$, $u_2$ thus defined are
  linearly independent solutions of (\ref{eq:severaleigen}). Further,
  their derivatives are given by uniformly convergent power series,
\begin{eqnarray}     
  u_1' &=& \frac{u_0'}{u_0}u_1 +  \frac{1}{pu_0} \sum_{\vn\ge0} 
     \frac{1}{(2|\vn|-1)!}  \sum_{i=1}^d \Xt{2\vn-\vd_i}
             \lambda_1^{n_1}  \cdots\lambda_2^{n_d} \nonumber  \\
  u_2' &=&  \frac{u_0'}{u_0}u_2 +  \frac{1}{pu_0}\sum_{\vn\ge0} 
              \frac{1}{(2|\vn|)!}\X{2\vn-\vd_i+\frac{1}{d}\vec1}
    \lambda_1^{n_1}\cdots\lambda_d^{n_i}  \label{eq:defu1'u2'} .
\end{eqnarray}
For every value of $\vec\lambda$, the initial values of these
functions are equal to
\begin{equation*}
  u_1(x_0) = u_0(x_0),\quad  u_1'(x_0) =  u_0'(x_0).\\
\end{equation*}
\begin{equation} \label{eq:initialcond}
  u_2(x_0) = 0,\quad   u_2'(x_0) = \frac{1}{p(x_0)u_0(x_0)}.
\end{equation} 
\end{theo}

\proof This proof is quite analogous to the proof for $d=1$ given in
\cite{KrP2010}, but certain details must be taken into account when
there are more spectral parameters. First we verify the convergence.

Let $\Lambda=\max(|\lambda_1|,\dots,|\lambda_d|)$. Recalling
(\ref{eq:M0Mi}), let $M=\max(M_0,M_1,\dots,M_d)$. Now by
(\ref{eq:Xtbound}),
\[ |\Xt{2|\vn|}(x)|\le \tilde c_{2|\vn|} M^{\left[|\vn|\right]} \cdot
  M^{\left[n_1\right]+\cdots+\left[n_d\right]} |x_2-x_1|^{2|\vn|}
 =  \tilde c_{\vj}  M^{|\vn|}|x_2-x_1|^{2|\vn|}
\]
so by (\ref{eq:ctildebound})  the summands in the formula for
$u_1$ in (\ref{eq:defu1u2}) are bounded by $a^{2|\vn|}/(2|\vn|)!$, where
\[  a=\sqrt{M\Lambda|x_2-x_1|}.
\]
Since $a^{2|\vn|}=a^{2n_1}a^{2n_2}\cdots a^{2n_d}$, we can factor the
sum $\sum_0^\infty a^{2|\vn|}/(2|\vn|)!$ into a product of $d$ sums, each of
which is equal to $\cosh a$. By comparison with this finite sum it
follows that the series for $u_1$ converges uniformly to a function
bounded by $\cosh^da$. By similar arguments the series for $u_1'$,
$u_2$, and $u_2'$ also converge uniformly, and this justifies the term
by term differentiation.

By Lemma \ref{lemm:Lu0Xt}, 
\begin{eqnarray*}
  Lu_1 &=& \sum_{|\vn|=0}^\infty  \frac{1}{(2|\vn|)!}
          \lambda_1^{n_1}\cdots\lambda_d^{n_d}  L(u_0\Xt{2\vn} ) \\
   &=&  \sum_{|\vn|=0}^\infty \frac{\lambda_1^{n_1}\cdots\lambda_d^{n_d}}{(2|\vn|)!} 
 (2|\vn|)(2|\vn|-1)u_0 \sum_{i=1}^d( r_i \Xt{2\vn-2\vd_i}  )    \\
   &=& u_0\sum_{|\vn|=0}^\infty \frac{\lambda_1^{n_1}\cdots\lambda_d^{n_d}}{(2|\vn|-2)!} 
   \sum_{i=1}^d r_i \Xt{2\vn-2\vd_i}     .
\end{eqnarray*}
Rearrange the last double sum as
\[  \sum_{i=1}^d r_i \sum_{|n|=0}^\infty 
     \frac{\lambda_1^{n_1}\cdots\lambda_d^{n_d}}{(2(|\vn|-1))!} 
     \Xt{2(\vn-\vd_i)}  
\]
and reindex each $n_i$ down by 1, using the assumption that $
\Xt{\vj}=0$ when $j_i<0$ for some $i$:
\[ \sum_{|\vn|=0}^\infty\frac{\lambda_1^{n_1}\cdots\lambda_d^{n_d}}{(2(|\vn|-1))!} 
     \Xt{2(\vn-\vd_i)}  =
    \sum_{|\vn|=0}^\infty
     \frac{\lambda_1^{n_1}\cdots\lambda_i^{n_i+1}\cdots\lambda_d^{n_d}}{(2(|\vn|))!} 
     \Xt{2\vn} 
\]\[ =
\lambda_i\sum_{|\vn|=0}^\infty\frac{\lambda_1^{n_1}\cdots\lambda_d^{n_d}}
   {(2(|\vn|))!}   \Xt{2\vn}   .
\]
  Thus we have
\[ Lu_1 = \left(\sum_i \lambda_i r_i \right) u_1 \] and the same
argument verifies the corresponding statement for $u_2$.  The final
statement regarding the initial values follows from the fact that
$\Xt{\vj}(x_0)=0$ whenever even a single $j_i$ is positive, so only
the constant terms survive in the series for $u_1$, $u_1'$; similarly
all but the lowest degree terms involving $\X{\vj}(x_0)$ also vanish.
\qed

It is well known that when $p,q,r_i$ are real-valued, a complex
nonvanishing solution $u_0$ of $Ly=0$ can be obtained as a complex
linear combination of any two linearly independent solutions; in fact,
by considerations of dimension one sees it is not necessary for the
coefficients be real-valued for such a nonvanishing solution to exist.
The hypotheses of Theorem \ref{theo:SPPS} could be weakened, for
instance by only requiring $1/(pu_0^2)$ and $r_iu_0^2$ to be
continuous, but we will not enter into such details here.

\begin{coro} \label{coro:SPPS} Let $u_1,u_2$ be given by (\ref{eq:defu1u2}).
Define
\begin{eqnarray*}
   v_1 &=& \frac{1}{u_0(x_0)}u_1 - p(x_0)u_0'(x_0)\,u_2, \nonumber\\
   v_2 &=& p(x_0)u_0(x_0)\,u_2.
\end{eqnarray*}
Then $v_1,v_2$ satisfy the normalizations
\begin{eqnarray*}
  v_1(x_0) = 1,\quad  &&v_1'(x_0) = 0, \nonumber\\
  v_2(x_0) = 0,\quad  && v_2'(x_0) = 1.
\end{eqnarray*}
\end{coro}
Observe that $v_1,v_2$ are also represented as power series in $\lambda_1,\dots,\lambda_d$.

\subsection{Generalized Sturm-Liouville equation}

We consider now equations of the form
\begin{equation}  \label{eq:slgeneralized}
  Ly = \sum_{i=1}^d \lambda_i R_i[y]
\end{equation}
where we define
\begin{equation}  \label{eq:Ri}
 R_i[y] = r_iy + s_i y'
\end{equation}
for given functions $r_i$, $s_i$, $i=1,\dots,d$. Thus
(\ref{eq:severaleigen}) is the particular case where all $s_i$ vanish
identically.  In \cite{KrT2015mod} SPPS formulas were developed for
(\ref{eq:slgeneralized}) for the case $d=1$.  The multiparameter
version is as follows. The formal powers are taken now with the same
starting values as previously but with the followng modified recursive
definition:
\begin{eqnarray}  
 \Xt{\vj} =  
 \DS |\vj| \int u_0R_i[u_0 \Xt{\vj-\vd_i}]  \  
      (\vj\mbox{ odd});  &&
 \Xt{\vj} =  
 \DS |\vj| \int \frac{1}{pu_0^2}\sum_{i=1}^d\Xt{\vj-\vd_i} \ 
      (\vj\mbox{ even}) \nonumber\\
 \X{\vj} =  
 \DS |\vj| \int u_0R_i[u_0 \X{\vj-\vd_i}]  \  
      (\vj\mbox{ even});  &&
 \X{\vj} =  
 \DS |\vj| \int \frac{1}{pu_0^2}\sum_{i=1}^d\X{\vj-\vd_i} \ 
      (\vj\mbox{ odd}); \label{eq:defXtgen}
\end{eqnarray} 
again we use integral entries in $\vj$ for $\Xt{\vj}$ and non-integral entries
for $\X{\vj}$.

We will need the common bound
\[ M = \sup_{[x_1,x_2]}\left(  \frac{1}{\left|pu_0^2\right|} ,\
       \big|u_0R_1[u_0]\big|,\ \dots,\   \big|u_0R_d[u_0]\big|,\ 
       \left|\frac{s_1}{p}\right|,\ \dots,\ \left|\frac{s_d}{p} \right| \right)
\]
of the functions appearing in the considerations below.  For odd $\vj$
it follows from  (\ref{eq:Ri})--(\ref{eq:defXtgen}) that
\begin{equation}  \label{eq:Riformula}  
   R_i[u_0\Xt{\vj-\vd_i}] =  R_i[u_0]\Xt{\vj-\vd_i} + 
     \frac{(|\vj|-1)s_i}{pu_0}\sum_{i'=1}^d\Xt{\vj-\vd_i-\vd_{i'}}
\end{equation}
which implies that the formal powers $\Xt{\vj}$ may be calculated without recourse to
numerical differentiation (other than for $u_0$) and that
\begin{equation}  \label{eq:Riestimate}  
 \left|u_0R_i[u_0\Xt{\vj-\vd_i}]\right| \le 
   M\bigg( |\Xt{\vj-\vd_i}| + (|\vj|-1)\sum_{i'=1}^d|\Xt{\vj-\vd_i-\vd_{i'}}| \bigg).
\end{equation}
Analogous statements hold for $\X{\vj}$.

\begin{lemm} \label{prop:generalized}  For all $x\in[x_2,x_2]$, the inequalities
  $|\Xt{\vj}|\le \widetilde P_{|\vj|}(x)$ and  $|\X{\vj}|\le P_{|\vj|}(x)$ hold, where
\begin{eqnarray*}
   \widetilde P_j(x) &=& d^{\left[\frac{j}{2}\right]} |\vj|! \!\!\! \sum_{k=[\frac{j}{2}]+1}^j
    {\left[\frac{j-1}{2}\right] \choose j-k} \frac{M^k}{k!}|x-x_0|^k  ,  \\
    P_j(x) &=& d^{\left[\frac{j-1}{2}\right]} |\vj|! \!\!\! \sum_{k=[\frac{j}{2}]-1}^j
    {\left[\frac{j}{2}\right] \choose j-k} \frac{M^k}{k!}|x-x_0|^k  ,
\end{eqnarray*}
for integral $j\ge0$. 
\end{lemm}
 
\proof First we consider $\Xt{\vj}$, i.e.\ $\vj$ has integer
entries. Write $E_k=(M|x-x_0|)^k/k!$ so $\left|M\int
  E_{k-1}\right|=E_k$.  The inequalities are clearly valid when $|\vj|$ is
0 or 1.  Suppose that it is valid for $|\vj|$ up to $n-1$. Now if $|\vj|=n$ is odd and
$\vj$ has an odd entry in the $i$-th position, we calculate that
\begin{eqnarray*}
  \widetilde P_{n-1}(x) &=&  d^{\frac{n-1}{2}} (n-1)!  \sum_{k=\frac{n+1}{2}}^{n-1}
    {\DS\frac{n-3}{2} \choose n-1-k } E_k ,\\
 d(n-1)| \widetilde P_{n-2}(x) )&=& d^{\frac{n-1}{2}}(n-1)!\sum_{k=\frac{n-1}{2}}^{n-2}
   {\DS\frac{n-3}{2} \choose n-2-k }  E_k  .
\end{eqnarray*}
Then by the inductive hypothesis and  (\ref{eq:defXtgen}), (\ref{eq:Riestimate}),
\begin{eqnarray*}
  \left|\Xt{\vj}(x)\right| &\le& n\int_{x_0}^x M(
    \widetilde P_{n-1}(x)+ d(n-1)  \widetilde P_{n-2}(x) ) \,dx  \\ 
 &=& d^{\frac{n-1}{2}} n!   \left( E_{\frac{n+1}{2}} + 
     \sum_{k=\frac{n+1}{2}}^{n-2} \bigg(
      {\DS\frac{n-3}{2} \choose n-1-k } + {\DS\frac{n-3}{2} \choose n-2-k }
        \bigg) E_{k+1}
      + E_{n} \right)  \\
 &=& d^{\left[\frac{n}{2}\right]}  n! \sum_{k=\frac{n-1}{2}}^{n-1} 
     {\DS\frac{n-1}{2} \choose n-1-k } E_{k+1}  \\
 &=& \widetilde P_{n}(x)  
\end{eqnarray*}
as is seen after reindexing  $k+1$ to $k$  and then noting that
$[n/2]+1=(n+1)/2$. On the other hand, if $\vj$ is even, then a
similar, simpler argument verifies the inequality.

 The verification for $\X{\vj}$ is analogous. \qed

 The following results for the generalized formal powers are now proved in
 exactly the same way as Lemmas \ref{lemm:Lu0Xt} and \ref{lemm:Lu0X}
 and Theorem \ref{theo:SPPS}.
\begin{lemm}
\[ L[u_0\Xt{2\vn}] = 2|\vn|(2|\vn|-1) \sum_{i=1}^d R_i[u_0\Xt{2\vn-2\vd_i}]
\]
and 
\[L(u_0\X{2\vn +\frac{1}{d}\vec1}) =   
     (2(|\vn|+1)(2|\vn|)\,  
    \sum_{i=1}^dR_i[u_0 \Xt{2\vn -2\vd_i+\frac{1}{d}\vec1}].
\]
\end{lemm}
\begin{theo} \label{theo:SPPSgen} Let
  $p,q,r_1,\dots,r_d,s_1,\dots,s_d$ be continuous on $[x_0,\,x_1]$,
  with $p$ continuously differentiable and $p(x)\not=0$. Define
  $\Xt{\vj}$ and $\X{\vj}$ by (\ref{eq:defXtgen}), and then define
  $u_1$, $u_2$ by (\ref{eq:defu1u2}). These series converge uniformly
  on $x\in[x_0,\,x_1]$ for every fixed
  $\vec\lambda=(\lambda_1\dots,\lambda_d)\in\C^d$, and are linearly
  independent solutions of the generalized Sturm-Liouville equation
  (\ref{eq:slgeneralized}). Their derivatives are given by
  (\ref{eq:defu1'u2'}) and they satisfy the initial conditions
  (\ref{eq:initialcond}).
\end{theo}

\subsection{Spectral problems}   

The treatment of multiparameter spectral problems by the SPPS approach
is the same as for a single spectral variable. Consider for simplicity
linear boundary conditions of the form
\begin{equation}  \label{eq:bdrycond}
  \alpha v(x_1) + \alpha' v'(x_1) = 0, \quad 
  \beta v(x_2) + \beta' v'(x_2) = 0 .
\end{equation}
For the general solution $v=c_1v_1+c_2v_2$ with $v_1$, $v_2$ given by
Corollary \ref{coro:SPPS}, this gives rise to a system of two
equations in $c_1,c_2$ with determinant
\[ \alpha(\beta v_1(x_2)+\beta v_1'(x_2) -
   \alpha'(\beta v_2(x_2)+\beta v_2'(x_2).
\]
Thus (\ref{eq:bdrycond}) is satisfied when $\chi(\vec\lambda)=0$, where
\begin{equation}  \label{eq:charfunc}
  \chi(\vec\lambda) = -\alpha'\beta\,v_1(x_2) + \alpha\beta\,v_1'(x_2)
  -\alpha'\beta'\,v_1'(x_2) + \alpha\beta'\,v_2'(x_2)   .
\end{equation}
Theorem \ref{theo:SPPS} represents $\chi(\vec\lambda)$ as a power
series in $\vec\lambda=(\lambda_1,\dots,\lambda_d)$. Solutions of the boundary
value problem are precisely the zeroes of this analytic function of
several complex variables.
   
 In like manner, nonlinear or mixed boundary conditions will also produce
a characteristic function.  When these conditions are analytic, the result
will be expressible as a power series in the $\lambda_i$, although it
may be more convenient to leave it as a function defined as a combination
of power series with other types of functions (cf.\ (\ref{eq:RT}) below).

Similarly, one may impose boundary conditions at more than two
points. One way of solving such a problem is by converting it
to an integral equation \cite{Ar1,Ar2,ChBo}.  With the approach
described here, one simply evaluates the SPPS representation at all
boundary points required, in order to obtain the desired set of
simultaneous characteristic equations.

\subsection{Remarks}

\subsubsection{Reduction to simple cases}
We note that for $d=1$ (i.e. $(\vec0)=(0)$, $(\vec1)=(1)$), the
starting integral of the $\Xt{\vj}$ family is $\Xt{0}=1$, and by
(\ref{eq:defX0}) for the family $\X{\vj}$ the starting integral also
reduces to
\[   \X{\frac{1}{1}\vec1-\vd_1} = \X{0} =\frac{1}{1} = 1.
\]  Further,  for general $d$ the degree-1 power $\X{(1/d)\vec1}$  is
simply the integral
\[    \int \frac{1}{pu_0^2}.
\]
which coincides with $\X{1}$ in the case $d=1$.  Thus our notation is
consistent with the ``classical'' definition of \cite{KrP2010}.

Considering $d>1$, let us suppose that $r_i$ is identically zero for every
$i\not=i_0$. Then $\Xt{\vj}$ will vanish whenever $\vj$ contains a
$j_i>0$ where $i\not=i_0$. The surviving powers
$\Xt{0,...0,,j_i,0,\dots,0}$ form the sequence $\Xt{j_i}$ of classical
1-spectral-parameter formal powers in the single variable
$\lambda_i$. Similarly, the $\X{\vj}$ reduce to the  sequence
$\X{j_i}$, and the series $u_1$, $u_2$ become the classical SPPS
solutions.
 
On the other hand, when all the $r_i$ are equal, the formal power
$\Xt{\vj}$ is unchanged when the entries $j_1,\dots,j_d$ are permuted,
so the sum only depends on the degree $|\vj|$, giving $\Xt{\vj}=
\tilde c_{\vj}\Xt{|\vj|}$, and similarly $\X{\vj}= c_{\vj}\X{|\vj|}$,
which are multiples of the classical formal powers.  It follows that
$u_1,u_2$ are the classical solutions obtained using
$|\vec\lambda|=\lambda_1+\cdots+\lambda_d$ in place of the single
spectral parameter.
 
\subsubsection{Computational aspects}

We make a few observations to simplify the task of programming the
formal powers.  One can omit the factors $|\vj|$ in the recursive
definitions (\ref{eq:defXt}), (\ref{eq:defX}) of $\Xt{\vj}$ and
$\X{\vj}$, producing ``rescaled powers'' $\stXt{\vj}$ and
$\stX{\vj}$ defined by
\[ \stXt{\vj} = \left\{ \begin{array}{ll}
 \DS  \int r_iu_0^2 \,\Xt{\vj-\vd_i},\quad  &
      \vj\mbox{ odd},  \\[2ex]
 \DS  \int \frac{1}{pu_0^2}\sum_{i=1}^d\Xt{\vj-\vd_i}, \quad&
      \vj\mbox{ even}.
 \end{array} \right.
\]
and similarly for $\stX{\vj}$. Then by induction
\[   \stXt{\vj} = \frac{1}{|\vj|!}  \Xt{\vj} ,\quad
  \stX{\vj} = \frac{1}{|\vj|!}  \X{\vj} .
\]
Besides this saving in multiplications when calculating the formal
powers (and often avoiding calculating with very large numbers), it
is no longer necessary to divide by these factorials to obtain the
terms in the sums for $u_1,u_2,u_1',u_2'$; i.e., we have simply
\[  u_1 = u_0\sum_{\vn}   \stXt{\vn} \lambda_1^{n_1}\cdots \lambda_d^{n_d} ,
\]
etc.  This is because the coefficent of each formal power
in the formulas (\ref{eq:defu1u2}) 
is precisely the reciprocal of the factorial of its degree.

The construction of the tables for $\tilde X$ and $X$ is seen to be
identical when we disregard the initial terms $\X{(1/d)\vec{1}-\vd_i}$
from the second table.  That is, according to whether we insert the
function $1=\Xt{\vec0}$ or $\int 1/(pu_0^2)=\X{(1/d)\vec1}$ in the
upper left hand corner, the same procedure of multiplying and then
integrating will produce the entire table for $\Xt{\vj}$ or $\X{\vj}$
respectively. Both tables and the corresponding power series can thus
be computed via a single program, except that in the formula (\ref{eq:defu1'u2'}) for 
$u_2'$, the first term corresponding to $\vn=\vec0$ contains negative
exponents and is not found in the truncated table.  Its value is
\[  \frac{1}{pu_0}\left( \frac{1}{0!}\sum_{i=1}^d\X{(1/d)\vec1-\vd_i}
\right) \lambda_1^0\cdots\lambda_d^0 =  \frac{1}{pu_0}.
\]
This term must be added in separately to obtain $u_2'$.

When programming, one may likely prefer to drop the fractional parts of
the indices, using effectively
\[  \ststX{\vj} =  \stX{\vj+(1/d)\vec1}.
\]
In the development of the theory given above, this amounts to
replacing the coefficient $|\vj|$ with $|\vj|+1$ , which is the
true degree of $\Xt{\vj}$.

It is easily seen that if $u_0$ is a solution of
(\ref{eq:severaleigen}) for a fixed multiparameter
$(\lambda_{1,0},\dots,\lambda_{d,0})\in\C^d$, then our construction of
$\Xt{\vj}$, $\X{\vj}$ will produce series in powers of
$\lambda_1-\lambda_{1,0}$, \dots, $\lambda_d-\lambda_{d,0}$ analogous
to (\ref{eq:defu1u2})--(\ref{eq:defu1'u2'}).  This can be used to
recenter the series for obtaining increased accuracy as in \cite{KrP2010}.

When calculating one must truncate the problem, say by using a finite
number $M$ of points of $[x_1,x_2]$ when integrating, and by
approximating the series (\ref{eq:defu1u2})--(\ref{eq:defu1'u2'}) with
polynomials formed of the terms for $|\vn|\le N$.  The total number of
formal powers in $\{\Xt{\vj},\X{\vj}\}_{|\vn|\le N}$ grows as $O(N^d)$, so the
memory requirement is of the order $O(MN^d)$.  For boundary value
problems this can be reduced by saving only the last value
$\Xt{\vj}(x_2)$, $\X{\vj}(x_2)$ once the values interior to the interval are
no longer needed for further integrations. The resulting memory cost
$O(MN)+O(N^d)$ is in fact a great savings since often $M$ is much
larger than $N$.

\section{Numerical examples}

We give some examples for $d=2$.  The operational parameters $M,N$ are
as described at the end of the last section; calculations were carried
out in \textit{Mathematica}.

\subsection{Boundary value problems}

\begin{figure}[!t]  \centering
 \pic{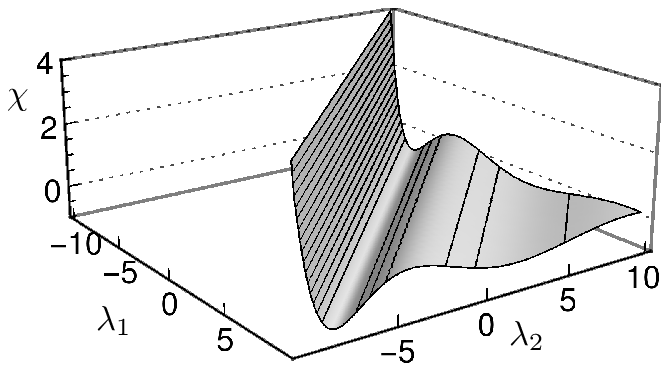}{-1,.1}{3.5}{}
 \pic{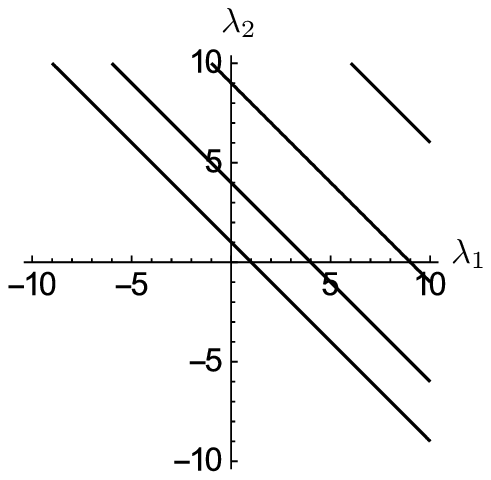}{6,0}{4}{scale=1}
  \caption{Characteristic function and zero-level curves (Example 1).}
  \label{fig:example1}
\end{figure}

\example This simple example uses constant coefficients $p=1$, $q=0$,
$r_1=r_2=-1$.  The equation $u''=-(\lambda_1+\lambda_2)u$ has
normalized solutions $v_1(x)=\cos(\sqrt{\lambda_1+\lambda_2}x)$,
$v_2(x)=\sin(\sqrt{\lambda_1+\lambda_2}x)/\sqrt{\lambda_1+\lambda_2}$.
On the interval $[x_1,x_2]=[0,\pi]$, the SPPS solutions of Corollary
\ref{coro:SPPS} with $M=800$, $N=20$ are found to agree with these
formulas to within $10^{-9}$ for $|\lambda_i|\le1$. As is common with
polynomial approximations, the accuracy drops rapidly for larger
values of $|\lambda_i|$ when the truncation limit $N$ is fixed.  We
impose the boundary conditions $u(0)=0$, $u(\pi)=0$.  The graph of the
characteristic function $\chi(\lambda_1,\lambda_2)$ (eigensurface) is
shown in Figure \ref{fig:example1}.  The eigencurves $\chi=0$,
calculated numerically from $\chi$ via the function
\texttt{ContourPlot} in the figure, coincide with the solutions of
\[   \lambda_1+\lambda_2 = \frac{k^2\pi^2}{b^2}
\] for $k=1,2,3,4$.
Indeed, the values of $|\chi(\lambda_1,\lambda_2)|$ for  $|\lambda_i|\le5$
for $k=1,2,3,4$ are less than $10^{-12}$,  $10^{-12}$,  $10^{-10}$,  $10^{-5}$
respectively.  When the maximal degree of the powers is reduced to $N=16$,
the level curve for $k=4$ is visibly far off the mark.

\begin{figure}[!b]  \centering
 \pic{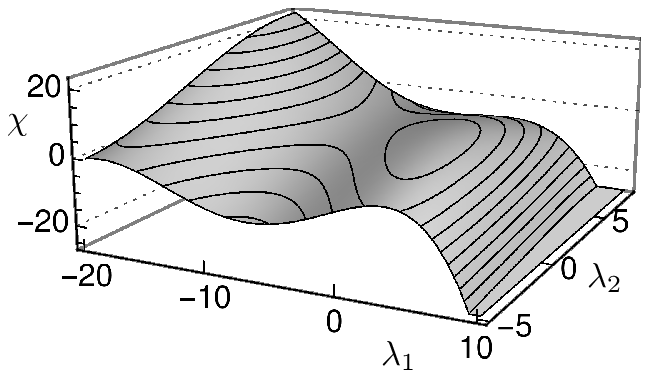}{-1,0}{3}{}
 \pic{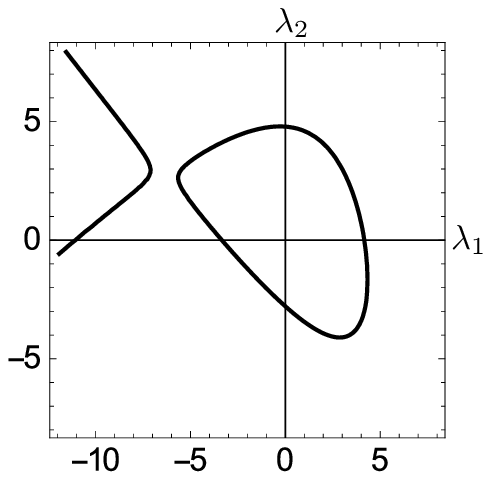}{6,0}{3}{scale=.9}
  \caption{Characteristic function and zero-level curves (Example 2).}
  \label{fig:example2}
\end{figure}

\begin{figure}[!t]  \centering
 \pic{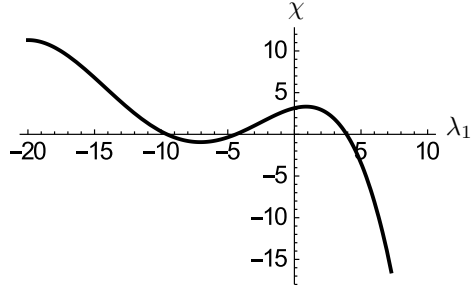}{2,0}{3}{}
   \caption{$\chi(\lambda_1,\lambda_2)$ for fixed value of $\lambda_2=1.0$.}
  \label{fig:example2b}
\end{figure}

\example This example, with $p(x)=1$, $q(x)=\cos x$,
$r_1(x)=\cos(x^2)$, $r_2(x)=\cos x$, which is not amenable to a
solution in closed form, is chosen to illustrate level sets which are
not connected and which contain closed curves. Using the same interval
$[0,\pi]$ and boundary conditions $u(0)=0$, $u(\pi)=0$, we find the
characteristic function and its zero set as depicted in Figure
\ref{fig:example2}.  For illustration we take an arbitrary section
$\lambda_2=1.0$, and restrict $\chi$ to this value (Figure
\ref{fig:example2b}). The corresponding numerical pairs
$(\lambda_1,\lambda_2)$ determine an ordinary differential equation
which can be solved numerically by \texttt{NDSolve} using the
boundary condition at $x=0$ to define an initial condition.  The resulting
values at $x=\pi$ were found to differ from
$\chi(\lambda_1,\lambda_2)$ by less than $10^{-6}$ when the experiment
was carried out with $M=100$, $N=12$.  The calculation of the
characteristic function took about $0.3$ seconds, and then each value
of $\lambda_1$ less than a thousandth of a second on a portable
computer (this does not include the time for checking by solving the
initial value problem).  The three eigenvalues
$\lambda_1\approx-9.5644,\ -4.3944,\ 3.9177$ in the range considered
are easily located by techniques of numerical approximation of zeroes
of polynomials.

  \begin{figure}[!b]  \centering
 \pic{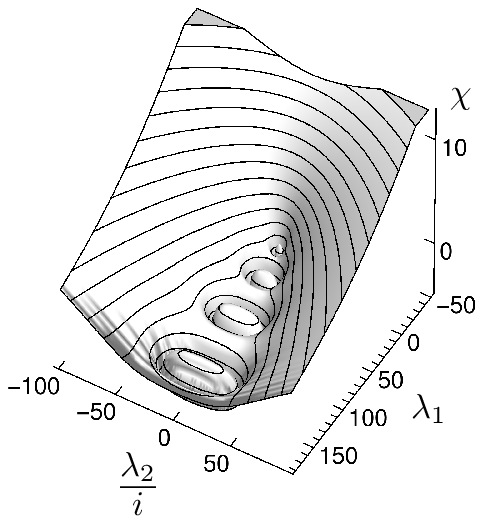}{-1,-.3}{4}{scale=1}
 \pic{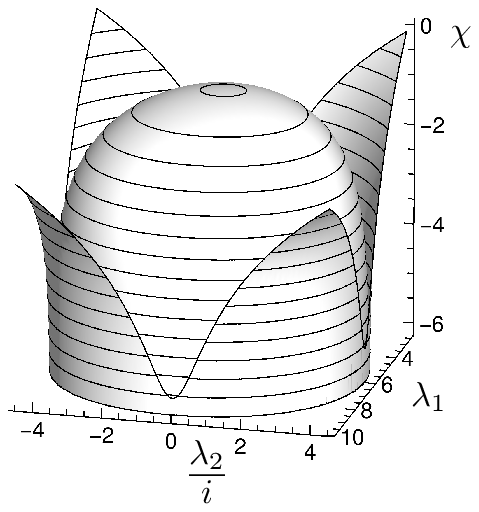}{5,-.3}{4}{scale=1}
 \caption{$\log|\chi(\lambda_1,\lambda_2)|$ for (\ref{eq:meissner})
   with $N=40$ (left); detail of region around smallest eigencurve
   with $N=16$ (right). }
  \label{fig:meissner}
\end{figure}

\example The following example involves consideration of complex eigenvalues.
The boundary value problem
\begin{equation}  \label{eq:meissner}
 y''(t) + (E + z\,\sgn t)y(t)=0,\quad y(-1)=y(1)=0,
\end{equation}
where $\sgn x$ is the sign of $x$, was studied in detail in \cite{MiVol}. A
spectral surface is formed of pairs $(E,z)\in\C^2$.  Let
$\lambda_1=E$, $\lambda_2=z$, $r_1(x)=-1$, $r_2(x)=\sgn(x)$ for $-1\le
x\le1$.

For the SPPS calculation, due to the jump singularity in $r_2$, it
would be appropriate to integrate separately on $[-1,0]$ and $[0,1]$.
For this example, however, we simply calculate the formal powers with
$M=10,000$ mesh points, and settle for about five significant figures
in the integrations.  Since $|\chi|$ does not change sign near its
zeros, we take the logarithm; then the plotting routine
(\texttt{Plot3D}) easily reveals the set where
$\chi(\lambda_1,\lambda_2)=0$ as shown in Figure \ref{fig:meissner}, where we
have taken $\lambda_1$ real and $\lambda_2$ purely imaginary.

In \cite{MiVol} certain curves in the spectral Riemann surface were
explicitly parametrized as
\begin{equation}  \label{eq:lam1lam2}
  \lambda_1(s)=s^2-h(s)^2,\quad \lambda_2(s)=2ish(s)
\end{equation}
where $s\in\bigcup_{n=0}^\infty[(n+1/2)\pi,(n+1)\pi]$ and where $h$
is defined implicitly by the relation
\[ s \sin(2s) + h(s)\sinh(2h(s)) =0.
\]
These curves were used to show that the surface is connected by
joining various 1-complex-dimensional parts. The first interval
$s\in[\pi/2,\pi]$ corresponds approximately to
$2.467\le\lambda_1\le\pi^2$, $0\le\lambda_2/i\le4.475$, and is the
smallest eigencurve revealed in the plot.  For the values given by
(\ref{eq:lam1lam2}) with values of $s$ in this interval we find
numerically that $|\chi(\lambda_1(s),\lambda_2(s))|<10^{-4}$ for
$s\in[\pi/2,\pi]$ when $N=20$.

\subsection{Application to electromagnetic transmission}

\example This example is based on \cite{CKKO} from
which we restate the minimum possible of background material.
The plane $\R^2=\{(x,y)\}$ is partitioned into the regions
\[  \Omega_1=\{ x<0\},\quad   \Omega_0=\{ 0<x<b\},\quad  \Omega_2=\{ x>b \}, 
\]
which are assumed to be composed of materials such that the index of
refraction in $\Omega_1$ and $\Omega_2$ takes constant values denoted
$n_1$, $n_2$ respectively, while in the inhomogeneous region
$\Omega_0$ it is a function $n=n(x)$ independently of $y$. These
values are bounded below by 1. An electromagnetic wave of the form
$e^{-ik_1x}$ travelling in $\Omega_1$ strikes the boundary line $x=0$
with $\Omega_0$ at an angle $\theta$ from the perpendicular, and is
partially reflected back into $\Omega_1$ as
$u(x)=e^{-ik_1x}+Re^{ik_1x}$ and partially transmitted into $\Omega_2$
at $x=d$ as $u(x)=e^{-ik_2x}$. The parameter
\begin{equation}  \label{eq:beta}
 \beta = k\sin\theta  
\end{equation} 
is introduced, where $k=2\pi/\lambda$ is the wave number in terms of
the wavelength $\lambda$ (here $\lambda$ will not denote an
eigenvalue).  In $\Omega_0$ the wave is governed by the differential
equation
\begin{equation}  \label{eq:optics}  
 u''(x)+(k^2n(x)^2-\beta^2)u(x)=0
\end{equation}
(for the ``s-polarization'', and a similar equation for the
``p-polarization''). The problem is the determination of the complex
constants $R$ and $T$, known as the reflection and transmission
coefficients.  In \cite{CKKO} the formulas
\begin{eqnarray}
  R &=& \frac{-k_1k_2v_2(b) - v_1'(b) - i k_2v_1(b)+ i k_1v_2'(b) }
           {(  v_1'(b)-k_1k_2v_2(b) ) + i(k_2v_1(b) + k_1v_2'(b)) }, \nonumber \\
  T &=& \frac{ 2ik_1(v_1(b)v_2'(b)-v_1'(b)v_2(b) )e^{-ik_2b} }
           {(  v_1'(b)-k_1k_2v_2(b) ) + i(k_2v_1(b) + k_1v_2'(b)) },  \label{eq:RT}
\end{eqnarray}
were given, where $k_1 = \sqrt{k^2n_1^2-\beta^2}$, $k_2 =
\sqrt{k^2n_2^2-\beta^2}$.  It was shown how by fixing $k$ in
(\ref{eq:beta}) and then using $\beta^2$ as the spectral parameter,
the SPPS formulas for dimension $d=1$ can be used to calculate $R$ and
$T$ for varying angles of incidence $\theta$.  Examples were given for
three sample functions $n(x)$.  All were for normal incidence
$\beta=0$, for which it is not difficult to calculate the solution of
the differential equation analytically in terms of special functions
for the examples considered (see for example \cite{Yeh}), and thus
compare the accuracy.  Similar calculations using SPPS were carried out in
\cite{Ced}, again for normal incidence, with many graphs comparing the
results to other numerical methods used in optics.

Equation (\ref{eq:severaleigen}) for $d=2$ with $\lambda_1=\beta^2$,
$\lambda_2=-k^2$, $r_1(x)=1$, $r_2(x)=n^2$ takes the form
(\ref{eq:optics}).  We apply Corollary \ref{coro:SPPS} to obtain
normalized solutions $v_1,v_2$, and then substitute these together
with
\begin{equation}
 k_1=\sqrt{-\lambda_1-\lambda_2 n_1^2},\quad
 k_2=\sqrt{-\lambda_1-\lambda_2 n_2^2}
\end{equation}
in (\ref{eq:RT}). This produces analytic expressions
$R(\lambda_1,\lambda_2)$, $T(\lambda_1,\lambda_2)$ which, while they
are not simple power series, serve conveniently for
calculations. 

\begin{figure}[!h]  \centering
 \pic{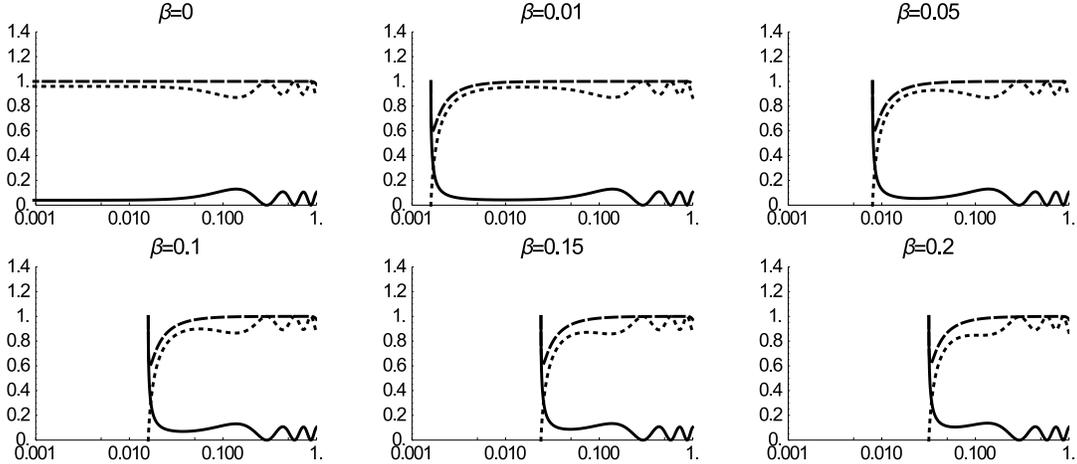}{-2,0.5}{6.5}{scale=.69}
 \caption{Logarithmic plots of $|R|^2$ (solid), $(n_2/n_1)|T|^2$
   (dotted) and formula (\ref{eq:rtcheck}) (dashed) as functions of
   the adimensional magnitude $b/\lambda_2\in[\beta/(2\pi),1]$, for
   $n(x)$ given by (\ref{eq:n(x)}).}
  \label{fig:optics}
\end{figure}
 
We will take one example, the ``hyperbolic'' refractive profile
\begin{equation}  \label{eq:n(x)}  
n(x) = n(0) e^{(x/b)\log(n(b)/n(0))} 
\end{equation}
with $d=1$, $n_1=1.0$, $n(0)=1.4$, $n(b)=2.1$, $n_2=1.5$. Further, we
set $b=1$. In Figure \ref{fig:optics} all graphs were plotted after a
single calculation of the series for $\chi(\lambda_1,\lambda_2)$ and
its substitution in the expressions (\ref{eq:RT}) for the parameters
given above.  It follows from (\ref{eq:beta}) that $b/\lambda\ge\beta
b/(2\pi)$, which determines our starting point for plotting the
curves.  For normal incidence $\beta=0$, conservation laws require the
expression
\begin{equation}  \label{eq:rtcheck} 
   |R|^2 + \frac{n_2}{n_1}|T|^2
\end{equation}
to be equal to 1; this is seen in the first graph, which agrees with
Figure 6 of \cite{CKKO}. For other values of $\beta$ we have
spot-checked numerically by selecting various values of the
dimensionless quantities $\beta$ and $b/\lambda$, then solving the
corresponding (\ref{eq:optics}) numerically with \texttt{NDSolve}, as
in the previous example.  The final values $v_1(b),\dots,v_2'(b)$
produce values of $R$, $T$ via (\ref{eq:RT}) for checking against the
$\chi$-values plotted here.  The results are given in Table
\ref{tab:opticserror}.  All of the data here is affected by the fact,
observed in \cite{CKKO}, that arithmetic operations in (\ref{eq:RT})
reduce the accuracy produced by the differential equations by several
significant figures.
 
\begin{table}[!htb]
  \centering
  \begin{tabular}{|l|c|c|c|c|c|} 
 \multicolumn{6}{c}{$M=30$, $N=10$}\\ \hline
 \parbox{7ex}{ \rule{0pt}{0pt}\hfill\raisebox{-4ex}{$\beta$}\\ $b/\lambda$  } 
          & 0.01    & 0.05    & 0.1    & 0.15    & 0.2    \\ \hline
 0.005    & 6       & 6       & 7      & 4       & 0      \\ 
 0.01     & ---     & 6       & 7      & 4       & 0      \\ 
 0.1      & ---     & ---     & 7      & 4       & 0      \\ 
 0.5      & ---     & ---     & 7      & 4       & 0      \\ 
 1        & ---     & ---     & 7      & 4       & 0      \\\hline  
  \end{tabular}
  \begin{tabular}{|l|c|c|c|c|c|}  
 \multicolumn{6}{c}{$M=50$, $N=16$}\\ \hline
 \parbox{7ex}{ \rule{0pt}{0pt}\hfill\raisebox{-4ex}{$\beta$}\\ $b/\lambda$  } 
          & 0.01    & 0.05    & 0.1    & 0.15    & 0.2    \\ \hline
 0.005    & 7       & 7       & 7      & 7       & 3      \\ 
 0.01     & ---     & 7       & 7      & 7       & 3      \\ 
 0.1      & ---     & ---     & 7      & 7       & 3      \\ 
 0.5      & ---     & ---     & 7      & 7       & 3      \\ 
 1        & ---     & ---     & 7      & 7       & 3      \\ 
 \hline
  \end{tabular}
  \caption{Number of significant digits in $\chi(\lambda_1,\lambda_2)$ 
    for selected values of $\beta$, $b/\lambda$.}
  \label{tab:opticserror}
\end{table}

 \section{Closing remarks}

We have shown how the representation of the solutions of the
Sturm-Liouville differential equation in terms of power series
in a single spectral parameter may be generalized to several
parameters $\lambda_1,\dots,\lambda_d$.  We hope that this will
make possible a deeper analysis and simplified computation
for many problems in physics and engineering, which have been
approached up to now by fixing the values of all parameters but
one, and solving by uniparameter methods.

Regarding the many aspects of uniparameter SPPS theory which have been
developed up to now, we point out as illustrative examples only two
possible areas for using several spectral parameters.

The so-called \emph{Sturm-Liouville pencils} 
\[   (pu')' + qu = (\sum r_i\lambda^i)u
\]
have been investigated from the SPPS perspective in
\texttt{arXiv:1401.1520}. This equation is a particular case of
(\ref{eq:severaleigen}) with $\lambda_1=1$, $\lambda_2=\lambda$,
\dots, $\lambda_d=\lambda^{d-1}$.  Thus our formulas provide the SPPS
series for this equation directly.

In another direction, coefficient functions with singularities at one of the
endpoints $[x_1,x_2]$, such as occur in Bessel's equation, have led to
modified versions of the SPPS formulas \cite{CKrT}.  Similar
results can be expected to hold also for several spectral parameters.
 
We close with the observation that an alternative construction to the
one described in this paper may be developed by first setting all but
one of the spectral parameters to zero, for example considering
\[ (py')' + qy = \lambda_1 r_1 y,
\]
and writing down the classical formulas for solutions $w_1^{[\lambda_1]}$,
$w_2^{[\lambda_1]}$ depending on this parameter.
These can be regarded as solutions of
\[  (py')' + (q - \lambda_1 r_1)y = 0,
\]
and after choosing a suitable nonvanishing linear combination, this
can be used as the seed for solving
\[  (py')' + (q - \lambda_1 r_1)y =  \lambda_2 r_2u
\]
to obtain $w_1^{[\lambda_1,\lambda_2]}$,
$w_2^{[\lambda_1,\lambda_2]}$, and so forth.  Even for the case $d=2$ the
resulting calculations to recover the coefficients of the SPPS series
turn out to be surprisingly complicated, and involve many products of
the nested integrals which cancel out at the end.  The author is
grateful to S.\ Torba for suggesting the simpler approach followed in
the present work.

This research was partially supported by grant 166183 of  
CONACyT (Consejo Nacional de Ciencia y Tecnolog{\'\i}a), Mexico.


\end{document}